\documentclass[11pt]{article}
\usepackage{amsmath,amssymb,amsthm, mathrsfs}

\begin{document}

\renewcommand{\baselinestretch}{1.5}
\renewcommand{\arraystretch}{1.5}

\begin{center}{\bf\LARGE Some Rules of Elementary Transformation  in $B^{+}(E,F)$ and their applications}\\
\vskip 0.5cm Ma Jipu$^{1,2}$
\end{center}

{\bf Abstract}\quad \small{Let $E,F$ be two Banach spaces, $B(E,F)$
 the set of all bounded linear operators from  $E$
into $F$, and $B^+(E,F)$ the set of  double splitting operators in
$B(E,F)$.
 In this paper, we present some rules  of elementary transformations  in $B^+(E,F)$,  consisting of five
 theorems. Let $\Phi_{m,n}$ be the set of all Fredholm operators $T$
 in $B(E,F)$ with dim$N(T)=m$ and codim$R(T)=n$, and $F_k=\{T\in
 B(E,F):$rank$T=k<\infty\}$. Applying the rules we  prove
 that  $F_k(k<$dim$F)$ and $\Phi_{m,n}(n>0)$ are path connected, so
 that they are not only smooth submanifolds in $B(E,F)$ with tangent
 space $M(X)=\{T\in B(E,F):TN(X)\subset R(X)\}$ at $X$ in them, but
 also path connected. Hereby we obtain the following
 topological construction: $B(\mathbf{R}^m,\mathbf{R}^n)=\bigcup^n\limits_{k=0}F_k(m\geq n),
B(\mathbf{R}^m,\mathbf{R}^n)=\bigcup^m\limits_{k=0}F_k(m<n),F_k$ is
path connected and smooth sub-hypersurface in
$B(\mathbf{R}^m,\mathbf{R}^n)(k<n),$ and  especially
dim$F_k=(m+n-k)k$ for $k=0,1,\cdots,n$. Finally we introduce an
equivalent relation in $B^+(E,F)$ and prove that the equivalent
class $\tilde{T}$ generated by $T$ is path connected for any $T\in
B^+(E,F)$ with $R(T)\varsubsetneq F$.

 }

{\bf Key words}\quad  Elementary Transformation Path Connected Set
of Operators Dimension Sub-hypersurface Smooth Submanifold
Equivalent Relation.


\vskip 0.2cm{\bf 1\quad Introduction }\vskip 0.2cm

Let  $E,F$ be two Banach spaces, $B(E,F)$  the set of all linear
bounded operators from  $E$ into $F$, and $B^+(E,F)$  the set of all
double splitting operators in $B(E,F)$.
 It is well known that the elementary transformation   of a matrix
 is a power tool in matrix theory.
 For example, by the elementary transformation  one can show that $F_k=\{T\in
 B(\mathbf{R}^n):$rank$T=k\}(k<n)$ is path connected. When
 dim$E=$dim$F=\infty$ there is no such  elementary transformation in
 $B(E,F)$, so it is difficult to prove that $F_k(k<$dim$F)$ in
 $B(E,F)$ is path connected. In this paper, we present some rules of
 elementary transformation in $B^+(E,F)$, which consist of five
 theorems, see  Section 2. Using the rules we prove that  $F_k(k<$dim$F)$ and
 $\Phi_{m,n}(n>0)$ are path connected,  where $\Phi_{m,n}$
 denotes the set of all Fredholm operators $T$ in $B(E,F)$ with
 dim$N(T)=m$ and codim$R(T)=n$. Then by Theorem 4.2 in [Ma4] we
 obtain the  following result  $F_k(k<$dim$F)$ and $\Phi_{m,n}(n>0)$
 are not only smooth submanifolds in $B(E,F)$ with tangent space
 $M(X)=\{T\in B(E,F):TN(X)\subset R(V)\}$ at any $X$ in then, but
 also path connected. As its application we have that
 $B(\mathbf{R}^m,\mathbf{R}^n)=\bigcup^n\limits_{k=0}F_k(m\geq n),
 B(\mathbf{R}^m,\mathbf{R}^n)=\bigcup^m\limits_{k=0}F_k(m<n),F_k(k<n)$,
 is path connected and smooth sub-hpersurface in
 $B(\mathbf{R}^m,\mathbf{R}^n)$ and especially, dim$F_k=(m+n-k)k$
 for $k=0,1,\cdots,n$. Finally we introduce an equivalent relation
 in $B^+(E,F)$ and prove that the equivalent class $\tilde{T}$
 generated by $T$ is path connected for $T\in B^+(E,F)$ with
 $R(T)\varsubsetneq F.$

 \vskip 0.2cm\begin{center}{\bf 2\quad Some Rules of Elementary Transformation  in $B^{+}(E,F)$}\end{center}\vskip 0.2cm


In this section, we  will introduce  some rules  of elementary
transformation  in $B^{+}(E,F)$, which consist of five theorems. It
is useful to imagene the trace of these rules of elementary
Transformations from Euclidean space to Banach space.

{\bf Theorem 2.1}\quad{\it  If $E = E_{1}\oplus R = E_{*}\oplus R$ , then the following conclusions hold:

$(i)$ there exists a unique $\alpha\in B(E_{*}, R)$ such that
$$
E_{1} = \left\{x + \alpha x: \forall x\in E_{*}\right\};\eqno(2.1)
$$}
conversely, for any $\alpha\in B(E_{*}, R)$  the subspace $E_{1}$ defined by $(2.1)$ satisfies $E = E_{1}\oplus R$

$(ii)$ $$
P^{R}_{E_{1}} = P^{R}_{E_{*}} + \alpha P^{R}_{E_{*}}\,\,\mathrm{and so}\,\,P_{R}^{E_{1}} = P_{R}^{E_{*}} - \alpha P^{R}_{E_{*}}.
$$

{\bf Proof}\quad For the proof of $(i)$ see [Ma3] and [Abr].

Obviously,

$\left(P^{R}_{E_{*}} + \alpha P^{R}_{E_{*}}\right)^{2} = P^{R}_{E_{*}} + \alpha P^{R}_{E_{*}},$ and
$$P^{R}_{E_{*}}x + \alpha P^{R}_{E_{*}}x = 0\,\, \mathrm{for}\,\, x\in E\Leftrightarrow P^{R}_{E_{*}}x=0
\Leftrightarrow x\in R$$
Then by (2.1), one concludes
$$
P^{R}_{E_{1}} = P^{R}_{E_{*}} + \alpha P^{R}_{E_{*}},\,\,\mathrm{and\, so,}\,\,P_{R}^{E_{1}} = P_{R}^{E_{*}} - \alpha P^{R}_{E_{*}}.
$$
The proof ends. \quad $\Box$

Let $B^{+}(E)$ be the set of all double splitting operators in $B(E)$ and $C_{r}(R) = \left\{T\in B^{+}(E): E = R(T) \oplus R\right\}$

{\bf Theorem 2.2}\quad{\it  Suppose that $E = E_{*}\oplus R$ and $\dim R > 0  .$  Then $P^{R}_{E_{*}}$ and $(-P^{R}_{E_{*}})$
are path connected in the set $\left\{T\in C_{r}(R): N(T) = R\right\}.$}

{\bf Proof}\quad Due to $\dim R > 0,$ one can assume that
$B(E_{*},R)$ contains a non-zero operator $\alpha$, otherwise
$E_*=\{0\}$ the theorem is trivial. Let $E_{1} = \left\{x+ \alpha
x:\forall x\in E_{*} \right\}.$ Then by Theorem 2.1, $E =
E_{1}\oplus R$ and
$$
P^{R}_{E_{1}} = P^{R}_{E_{*}} + \alpha P^{R}_{E_{*}}. \eqno(2.2)
$$
Consider the path
$$
P(\lambda)= (1-2\lambda)
P^{R}_{E_{*}}+(1-\lambda)P^{R}_{E_{*}}\,\,\,\, 0\leq\lambda\leq1.
$$
Clearly,
$$
R(P(\lambda))= R(P^{R}_{E_{*}} + \frac{1-\lambda}{1-2\lambda}\alpha
P^{R}_{E_{*}}),\quad\,\,\,0\leq\lambda\leq1.
$$
Then by Theorem 2.1, $R(P(\lambda))\oplus R = E$, i.e.$P(\lambda)\in
C_{r}(R)$ for all $\lambda\in [0,1]$. In addition, $P(0)=P^R_{E_1},
P(1)=-P^R_{E_*}$ and $N(P(\lambda))=R,\forall\lambda\in[0,1]$. This
shows that $P^{R}_{E_{1}}$ and $-P^{R}_{E_{*}}$ are path connected
in $\left\{T\in C_{r}(R): N(T) = R\right\}.$

Next go to show that $P^{R}_{E_{1}}$ and $P^{R}_{E_{*}}$ are path connected in
$\left\{T\in C_{r}(R): N(T) = R\right\}.$ Consider the path
$$
P(\lambda) =  P^{R}_{E_{*}} + \lambda\alpha P^{R}_{E_{*}}\,\,\,\,0\leq\lambda\leq1.
$$
By Theorem 2.1, $P(1) =  P^{R}_{E_{*}} + \alpha P^{R}_{E_{*}}=  P^{R}_{E_{1}},$
 and $R(P(\lambda))\oplus R = E\,\,\,\,\forall \lambda\in [0,1],$ where $E_{1} = \left\{x+\alpha x: \forall x\in E_{*}\right\}.$
Obviously, $N(P(\lambda))= R$ and $P(0)= P^{R}_{E_{*}}$. Therefor
$P^{R}_{E_{1}}$ and $P^{R}_{E_{*}}$ are path connected in
$\left\{T\in C_{r}(R): N(T) = R\right\}.$ Thus the theorem is
proved. \quad $\Box$

For simplicity, still write $C_{r}(N)=\left\{T\in B^{+}(E,F):
R(T)\oplus N= F\right\}$ in the sequal.

{\bf Theorem 2.3}\quad{\it  Suppose  $T_{0}\in C_{r}(N)$  and $F =
F_{*}\oplus N$. Then $T_{0}$ and $P^{N}_{F_{*}}T_{0}$ are path
connected in the set $\left\{T\in C_{r}(N): N(T) =
N(T_{0})\right\}.$}

{\bf Proof}\quad One can assume $R(T_{0}) \neq F_{*}$, otherwise the
theorem is trivial. Then by Theorem 2.1, there exists a non-zero
operator $\alpha\in B(F_{*},N)$ such that
$$
R(T_{0})=\left\{y+\alpha y : \forall y\in F_{*}\right\},\,\,\,P^{N}_{R(T_{0})}= P^{N}_{F_{*}}+\alpha P^{N}_{F_{*}},
$$
and
$P_{N}^{R(T_{0})}= P_{N}^{F_{*}}-\alpha P^{N}_{F_{*}}.$ So
$$
T_{0}=\left\{P^{N}_{F_{*}}+\alpha P^{N}_{F_{*}}\right\}T_{0}.\eqno(2.3)
$$
Let
$$
F_{\lambda}=\left\{y+\lambda\alpha y : \forall y\in
F_{*}\right\}\,\,\,\,\mathrm{for\ all}\,\,\,\,\lambda\in [0,1].
$$
Note  $\lambda\alpha\in B(F_{*},N).$ According to Theorem 2.1 we
also have
$$
P^{N}_{F_{\lambda}}=P^{N}_{F_{*}}+\lambda\alpha
P^{N}_{F_{*}}\quad\mbox{and}\quad
P^{F_{\lambda}}_{N}=P^{F_{*}}_{N}-\lambda\alpha P^{N}_{F_{*}}
,\quad\quad \forall\lambda\in[0,1].$$

Consider the path
$$
P(\lambda)=P^{N}_{F_{\lambda}}T_{0},\;\;\forall \lambda\in[0,1].
$$

Because $F=R(T_{0})\oplus N=F_{\lambda}+N$ one observers
$$
R(P(\lambda))=F_{\lambda},\;\;\forall \lambda\in[0,1],
$$
and so
$$
R(P(\lambda))\oplus N=F,\quad{\rm i.\ e.,}\quad  P(\lambda)\in
C_r(N).
$$
Note
$$
y\in N(P(\lambda))\Leftrightarrow P^{N}_{F_{\lambda}}T_{0}y=0\Leftrightarrow T_{0}y\in N\Leftrightarrow y\in N(T_{0})
$$
i.e., $N(P(\lambda))=N(T_{0}), \forall \lambda\in[0,1].$ Thus
$$
P(\lambda)\in \{T\in C_{r}(N):N(T)=N(T_{0})\},\;\lambda\in [0,1].
$$
In addition, $P(1)=T_{0}$ by (2.3), and $P(0)=P^{N}_{F_{*}}T_{0}$.

Finally we conclude that $T_{0}$ and $P^{N}_{F_{*}}T_{0}$ are path
connected in the set  $\{T\in C_{r}(N):N(T)=N(T_{0})\}$. The proof
ends.\quad$\Box$


Let $C_{d}(R)=\{T\in B^{+}(E,F):E=N(T)\oplus R\}$.

{\bf Theorem 2.4}\quad {\it Suppose that $T_{0}\in C_{d}(R_{0})$ and
$E=E_{*}\oplus R_{0}.$ Then $T_{0}$ and $T_0P^{E_{*}}_{R_{0}}$ are
path connected in the set $\{T\in C_{d}(R_{0}):R(T)=R(T_{0})\}$}.

\textbf{Proof}\quad  One can assume $E_{*}\neq N(T_{0}),$ otherwise
the theorem is trivial. Then by Theorem 2.1, there exists a non-zero
operator $\alpha\in B(E_{*},R_{0})$ such that
$$
P^{R_{0}}_{N(T_{0})}=P^{R_{0}}_{E_{*}}+\alpha P^{R_{0}}_{E_{*}},
$$
and so,
$$
T_{0}=T_{0}P^{N(T_{0})}_{R_{0}}=T_{0}(P^{E_{*}}_{R_{0}}-\alpha
P^{R_{0}}_{E_{*}}).\eqno(2.4)
$$

Consider the path as follows,
$$
P(\lambda)=T_{0}(P^{E_{*}}_{R_{0}}-\lambda\alpha
P_{E_{*}}^{R_{0}}),\;\;0\leq\lambda\leq1.
$$
Since $(P^{E_{*}}_{R_{0}}-\lambda\alpha
P_{E_{*}}^{R_{0}})x=x\;\forall x \in R_{0},$ we conclude
$R(P(\lambda))=R(T_{0})$. We also have
$N(P(\lambda))=R(P^{E_{*}}_{R_{0}})+\lambda\alpha
P^{E_{*}}_{R_{0}})$.Indeed,
$$
x\in N(P(\lambda))
\Longleftrightarrow (P^{E_{*}}_{R_{0}}-\lambda\alpha P^{R_{0}}_{E_{*}})x=0
\Longleftrightarrow x\in R(P^{R_{0}}_{E_{*}}+\lambda\alpha P^{R_{0}}_{E_{*}})
$$
since  $\lambda\alpha\in B(E_{*},R_{0})$ for all $\lambda\in[0,1].$

Thus  $P(\lambda)\in \{T\in C_{d}(R):R(T)=R_{0}\},
\forall\lambda\in[0,1].$ In addition, $P(1)=T_{0}$ by (2.4), and
$P(0)=T_{0}P^{E_{*}}_{R_{0}}.$ Then the theorem is
proved.\quad$\Box$


{\bf Theorem 2.5}\quad{\it Suppose that the subspaces $E_{1}$ and
$E_{2}$ in $E$ satisfy $dim E_{1}=dim E_{2}< \infty.$ Then $E_{1}$
and $E_{2}$  possess a common complement R, i.e.,$E=E_{1}\oplus R=
E_{2}\oplus R$.}

\textbf{Proof}\quad  According to the assumption $dim E_{1}=dim
E_{2}<\infty$, we have the following decompositions:
$$
E=H\oplus(E_{1}+E_{2}),\;E_{1}=E^{*}_{1}\oplus(E_{1}\cap E_{2}),\;\mbox{and}\;E_{2}=E^{*}_{2}\oplus(E_{1}\cap E_{2}).
$$
It is easy to observe that  $(E^{*}_{1}\oplus
E^{*}_{2})\cap(E_{1}\cap E_{2})=\{0\}$ and $dim E^{*}_{1}=dim
E^{*}_{2}<\infty.$ Indeed, if $e^{*}_{1}+e^{*}_{2}$ belongs to
$E_{1}\cap E_{2}$, for $e^{*}_{i}\in E^{*}_{i},\;i=1,2,$
 then
$e_{2}^{*}=(e^{*}_{1}+e^{*}_{2})-e^{*}_{1}\in E_{1}$ and
$e_{1}^{*}=(e^{*}_{1}+e^{*}_{2})-e^{*}_{2}\in E_{2}$, so that
$e^{*}_{1}=e^{*}_{2}=0$ because of $E^{*}_{i}\cap(E_{1}\cap
E_{2})=\{0\},\,i=1,2.$ Hereby, one can see
$$
E_{1}+ E_{2}=(E^{*}_{1}\oplus E^{*}_{2})\oplus(E_{1}\cap
E_{2}).\eqno(2.5)
$$
we now are in the position to determine $R$. We may assume $dim
E^{*}_{1}=dim E^{*}_{2}>0,$ otherwise the theorem is trivial. Then
$B^{X}(E^{*}_{1},E^{*}_{2})$ contains an operator $\alpha$, which
bears the subspace $H_{1}$ as follows,
$$
H_{1}=\{x+\alpha x:\forall x\in E^{*}_{1}\}=\{x+\alpha^{-1}x:\forall x\in E^{*}_{2}\}.
$$
By Theorem 2.1,
$$
E^{*}_{1}\oplus E^{*}_{2}=H_{1}\oplus E^{*}_{2}=H_{1}\oplus E^{*}_{1}
$$
Finally, according to (2.5) we have
$$
E=H\oplus(E_{1}+ E_{2})=H\oplus(E^{*}_{1}\oplus E^{*}_{2})\oplus(E_{1}\cap E_{2})=H\oplus H_{1}\oplus E_{2}^{*}\oplus(E_{1}\cap E_{2})=H\oplus H_{1}\oplus E_{2}
$$
and
$$
E=H\oplus H_{1}\oplus E_{1}^{*}\oplus(E_{1}\cap E_{2})=H\oplus H_{1}\oplus E
$$
This says $R=H\oplus H_{1}.$ The proof ends.\quad$\Box$


\vskip 0.2cm\begin{center}{\bf 3\quad Some
Applications}\end{center}\vskip 0.2cm

In this section we  will give some application of  the rules in
Section 2.


{\bf Theorem 3.1}\quad{\it For $k < dim F$, $F_k$ is  path
connected.}

{\bf Proof}\quad  In what follows, we may assume $k>0$, otherwise,
the theorem is trivial. Let $T_1$ and $T_2$ be arbitrary two
operators in $F_k$, and
$$E=N(T_{i})\oplus R_{i},\quad\quad\quad i=1,2.$$
Then $0<$ dim $R_{1}=$dim $R_{2}=k<\infty$, and  by Theorem 2.5,
there exists a subspace $N_{0}$
 in $E$ such that
$$E=N_{0}\oplus R_{1}=N_{0}\oplus R_{2}\eqno(3.1)$$
Let
$$L_{i}x=\left\{\begin{array}{rcl}T_{i}x, & & {x\in R_{i}}\\
0, & & {x\in N_{0}}
\end{array}\right.$$
$i=1,2.$

We claim that $T_{i}$ and $L_{i}$ are path connected in $F_{k},\,i=1,2.$ Due to $(3.1)$ the theorem 2.1 shows that there exists an operator
$\alpha_{i}\in B(N_{0},R_{i})$ such that
$$
P^{N(T_{i})}_{R_{i}} = P^{N_{0}}_{R_{i}}- \alpha_{i}P^{R_{i}}_{N_{0}},\; i=1,2.
$$

Consider the following paths:
$$ P_{i}(\lambda)= T_{i}(P^{N_{0}}_{R_{i}}- \lambda\alpha_{i}P^{R_{i}}_{N_{0}}),\;0\leq \lambda \leq 1 \,\mbox{and}\, i = 1,2.$$
Obviously, $P_{i}(0)=T_{i}P^{N_{0}}_{R_{i}}=L_{i},$ and
$P(1)=T_{i}(P^{N_{0}}_{R_{i}}-\alpha_{i}P^{R_{i}}_{N_{0}})=
T_{i}P^{N(T_{i})}_{R_{i}}=T_{i},\;i=1,2.$ Next go to show
$R(P_i(\lambda))=T_{i}R_{i}=R(T_{i}).$

Evidently, $R(P^{N_{0}}_{R_{i}}-
\lambda\alpha_{i}P^{R_{i}}_{N_{0}})\subset R_{i},$ and the converse
rlation follows from $(P^{N_{0}}_{R_{i}}-
\lambda\alpha_{i}P^{R_{i}}_{N_{0}})r = r, \forall r\in R_{i}.$ This
says $P_{i}(\lambda)\in F_{k}$ for $0\leq \lambda \leq 1
\,\mbox{and}\, i = 1,2,$ so that $L_{i}$ and $T_{i}$ are path
connected in $F_{k},\; i=1,2.$ Hence, in what follows, we can assume
$N(T_{1})=N(T_{2})=N_{0}.$ In the  other hand, according to Theorem
2.5 we have
$$F = R(T_{1})\oplus N_{*}= R(T_{2})\oplus N_{*}\eqno(3.2)$$
where $ N_{*}$ is a closed subspace in $F$. Then by Theorem 2.3, the
proof of the theorem turns to that of the following conclusion:
$P^{N^{*}}_{R(T_{1})}T_{2} $ and $T_{1}$ are path connected in
$F_{k}$. Let
$$T^+_{1}y=\left\{\begin{array}{rcl}\left(T_{1}|_{R_{1}}\right)^{-1}y, & & {y\in R(T_{1})}\\
0, & & {y\in N_{*}}
\end{array}\right.$$
It is easy to observe
$$
P^{N_{*}}_{R_{T_{1}}}T_{2} =
\left(P^{N_{*}}_{R_{(T_{1})}}T_{2}T^{+}_{1}\right)T_{1}.\eqno(3.3)
$$
In fact, note $N(T_{1})=N(T_{2})=N_{0}$ and $T^{+}_{1}T_{1}=P^{N_{0}}_{R_{1}},$ then
$$
 P^{N_{*}}_{R_{(T_{1})}}T_{2} = P^{N_{*}}_{R_{(T_{1})}}T_{2}\left(P^{N_{0}}_{R_{1}}+P_{N_{0}}^{R_{1}}\right)
  = P^{N_{*}}_{R_{(T_{1})}}T_{2}P^{N_{0}}_{R_{(T_{1})}}
  = P^{N_{*}}_{R_{(T_{1})}}T_{2}T^{+}_{1}T_{1}.
$$

We  claim $P^{N_{*}}_{R_{(T_{1})}}T_{2}T^{+}_{1}|_{R(T_{1})}\in
B^{X}(R(T_{1})).$

Evidently,
\begin{eqnarray*}
&& P^{N_{*}}_{R_{(T_{1})}}T_{2}T^{+}_{1}y=0\,\, \mbox{for}\, y\in R(T_{1})
 \Leftrightarrow  T_{2}T^{+}_{1}y\in N_{*}\\
&& \Leftrightarrow  T^{+}_{1}y\in N_{0}\cap R_{1} \,\,(\mbox{because of (3.2)})\\
&& \Leftrightarrow  T^{+}_{1}y=0\Leftrightarrow y=0,
\end{eqnarray*}
i.e.,$P^{N_{*}}_{R_{(T_{1})}}T_{2}T^{+}_{1}|_{R(T_{1})}$ is injective.

In the  other hand, since the decomposition (3.2) implies
$P^{N_{*}}_{R_{(T_{1})}}|_{R(T_{2})}\in B^{X}\left(R(T_{2}),
R(T_{1})\right)$, there is a unique $r_{2}\in R_{2}$ for any  $y\in
R(T_{1})$, such that
$$
y
= P^{N_{*}}_{R_{(T_{1})}}T_{2}r_{2}
= P^{N_{*}}_{R_{(T_{1})}}T_{2}\left(P^{N_{0}}_{R_{1}}r_{2}+P_{N_{0}}^{R_{2}}r_{2}\right)
= P^{N_{*}}_{R_{(T_{1})}}T_{2}P^{N}_{R_{1}}r_{2}
$$
Then $y_{0} = T_{1}P^{N_{0}}_{R_{1}}r_{2}$ fulfills
$$
y = P^{N_{*}}_{R_{(T_{1})}}T_{2}P^{N_{*}}_{R_{1}}r_{2} = P^{N_{*}}_{R_{(T_{1})}}T_{2}T_{1}^{+}y_{0}.
$$
This says that $P^{N_*}_{R(T_1)}T_2T^{+}_1|_{R(T_1)}$ is surjective.
So $P^{N_*}_{R(T_1)}T_2T^{+}_1|_{R(T_1)}\in B^X(R(T_1))$ of all
invertible  operators  in $B(R(T_1)$.

It is well known that $P^{N_*}_{R(T_1)}T_2T^{+}_1|_{R(T_1)}$  is
path connected with some one of $-I_{R(T_1)}$ and $I_{R(T_1)}$ is
$B^X(R(T_1))$,  where $I_{R(T_1)}$ denotes  the identity on
$R(T_1)$.
 Let $Q(t)$ be a path in $B^X(R(T_1))$ with
$Q(0)=P^{N_*}_{R(T_1)}T_2T_1^+|_{R(T_1)}$ and $Q(1)=-I_{R(T_1)}$ (or
$I_{R(T_1)})$. Then $Q_1(t)=Q(t)T_1$ is a path in the set $S=\{T\in
B(E,F):R(T)=R(T_1)$ and $N(T)=N(T_1)\}$ satisfying
$$Q_1(0)=P^{N_*}_{R(T_1)}T_2T^+_1T_1\quad{\rm and}\quad
Q_1(0)=-T_1({\rm or}\ T_1).$$ Since $S\subset F_k$, it follows that
$P^{N_*}_{R(T_1)}T_2T^+_1T_1$ and $-T_1($or $T_1)$ are path
connected in $F_k$. By (3.3) we merely need to show the following
conclusion:  if $P^{N_*}_{R(T_1)}T_2T^+_1T_1$ is path connected with
$-T_1$ in $F_2$ then  $P^{N_*}_{R(T_1)}T_2T^+_1T_1$ and $T_1$ are
path connected in $F_k$. Note $F=R(T_1)\oplus N_*$ and dim$N_*>0$.
By Theorem 2.2, $P^{N_*}_{R(T_1)}$ and $-P^{N_*}_{R(T_1)}$ are path
connected in the set $\{T\in G_r(N_*):N(T)=N_*\}=\{T\in
B(F):F=R(T)\oplus N_*$ and $N(T)=N_*\}$. Let $Q(t)$ be such a  path
in the set with $Q(0)=-^{N_*}_{R(T_1)}$ and $Q(1)=P^{N_*}_{R(T_1)}.$
Since $F=R(Q(t))\oplus N_*$ and $N(Q(t))=N_*\forall t\in[0,1],
R(Q(t)T_1)=Q(t)R(T_1)=Q(t)F,$  $Q(t)$ for $t\in[0,1]$ belongs to
$F_k$, and satisfies $Q(0)T_1=-P^{N_*}_{R(T_1)}T_1=-T_1$ and
$Q(1)T_1=T_1.$ This says that $-T_1$ and $T_1$ are path connected in
$F_k$, so that the conclusion is proved. \quad$\Box$


It is obvious that if either dim$E$ or dim$F$ is finite, then
$B(E,F)$ consists of all finite rank operators. Hence we assume
dim$E=$dim$F=\infty$ in the sequel.

 Let $\Phi_{m,n} = \left\{T\in
B^{+}(E,F): \dim N(T) = m <\infty\; \mbox{and} \;\mathrm{codim} R(T)
= n < \infty\right\}$, we have


{\bf Theorem 3.2}\quad{\it  $\Phi_{m,n} , (n>0)$ is  path
connected.}

{\bf Proof}\quad Let $T_{1}$ and $T_{2}$ be arbitrary two operators
in $\Phi_{m,n}$ and
$$
F = R(T_{1}) \oplus N_{1} = R(T_{2}) \oplus N_{2},
$$
i.e. $T_{1}\in C_{r}(N_{1})$ and $T_{2}\in C_{r}(N_{2})$.

Clearly, $\dim N_{1} = \dim N_{2} = n<\infty$. Then by Theorem 2.5,
there exists a subspace $F_{*}$ in $F$ such that
$$
F = R(T_{1}) \oplus N_{1} =F_{*} \oplus N_{1}\, \mathrm{and}\, F = R(T_{2}) \oplus N_{2} =F_{*} \oplus N_{2}.\eqno(3.4)
$$
So, by Theorem 2.3 the proof of the theorem turns to that of the
operators $P^{N_{1}}_{F_{*}}T_{1}$ and $P^{N_{2}}_{F_{2}}T_{2}$
being path connected in $\Phi_{m,n}.$ For simplicity, still write
$P^{N_1}_{F_1}T_1$ and $P^{N_2}_{F_*}T_2$ as $T_1$ and $T_2$,
respectively. However, here $R(T_1)=R(T_2)=F_*$ while  $N(T_1),
N(T_2)$ keep invariant. Due to dim$N(T_1)=$dim$N(T_2)=m<\infty.$
Theorem 2.5 shows that there exists a subspace $R$ in $E$ such that
$$E=N(T_1)\oplus R=N(T_2)\oplus R.\eqno(3.5)$$
Then by Theorem 2.4, one can conclude that $T_2P^{N(T_1)}_R$ and
$T_2$ are path connected in $\Phi_{m,n}.$ Thus the proof of the
theorem turns once more to that of $T_2P^{N(T_1)}_R$ and $T_1$ being
path connected in $\Phi_{m,n}.$ In what follows, we do this.

Let
$$T_{1}^{+}y =\left\{\begin{array}{rcl}\left(T_{1}|_{R}\right)^{-1}y, & & {y\in F_{*}},\\
0, & & {y\in N_{1}}.
\end{array}\right.$$
Then
$$
T_{2}P^{N(T_{1})}_{R} = T_{2}T^{+}_{1}T_{1}.\eqno(3.6)
$$
It is easy to observe $T_{2}T^{+}_{1}\mid _{F_{*}}\in B^{X}(F_{*}).$ Indeed,
\begin{eqnarray*}
&& T_{2}T^{+}_{1}y=0\,\, \mbox{for}\, y\in F_{*}
 \Leftrightarrow  T^{+}_{1}y\in N(T_{2})\cap R\\
&& \Leftrightarrow  T^{+}_{1}y = 0 \,\,\mbox{because of (3.5)} \Leftrightarrow   y=0,
\end{eqnarray*}
i.e., $T_{2}T^{+}_{1}|_{F_{*}}$ is injective; let;  $r\in R$ satisfy
$T_2r=y$ for any $y\in F_*$, and $T_1r-y_0\in F_*$, then
$$T_2T^+_1y_0=T_2T^+_1T_1r=T_2P^{N(T_1)}_Rr=T_2r=y,$$
and so $T_2T^+_1|_{F_*}$ is also surjective. Thus
$T_2T^+_1|_{F_*}\in B^X(F_*)$ and is path connected with some one of
$-I_{F_*}$ and $I_{F_*}$ is in $B^X(F_*)$. Similar to the way in the
proof of Theorem 3.1 one can prove that $T_2T^+_1T_1$ is path
connected with some one of $-T_1$ and $T_2$ is in $\Phi_{m,n}$. Note
the  equality $F=F_*\oplus N_1$ in (3.4) and dim$N_1>0.$ Due to dim
$N_1>0$ Theorem 2.2 shows that $-P^{N_1}_{F_*}$ and $P^{N_1}_{F_*}$
are path connected in the set $S=\{T\in C_F(N_1)\subset (F);
N(T)=N_1\}$, say that $Q(X)\in S$ for $\lambda\in[0,1]$ fulfills
$Q(0)=-P^{N_1}_{F_*}$ and $Q(1)=P^{N_1}_{F_*}.$ Obviously,
$$R(Q(\lambda)T_1)=Q(\lambda)F_*=Q(\lambda)(F_*\oplus
N_1)=R(Q(\lambda))\quad({\rm note}\  R(T_1)=F_*),$$ and
$$Q(\lambda)T_1x=0\Leftrightarrow T_1x\in N_1\Leftrightarrow x\in
N(T_1)\quad({\rm note}\ F_*\cap N_1=\{0\})$$ for all
$\lambda\in[0,1].$ So $Q(\lambda)T_1$ for any $\lambda\in[0,1]$
belongs to $\Phi_{m,n}$. This says that $-T_1$ and $T_1$ are path
connected in $\Phi_{m,n}$. Therefore by (3.6) $T_2P^{N(T_1)}_R$ is
also  path connected with $T_1$ in $\Phi_{m,n}$.$\Box$

By Theorem 4.2 in [Ma4] we futher have

\textbf{Theorem 3.3}\quad $F_k(k<$dim$F$) and $\Phi_{m,n}(n>0)$ are
not only path connected set  but also smooth submanifolds in
$B(E,F)$ with tangent space $M(X)=\{T\in B(E,T):TN(X)\subset R(X)\}$
at any $X$ in them.

Applying the theorem to $B(\mathbf{R}^m,\mathbf{R}^n)$ we obtain the
geomatrical and topologial construction of
$B(\mathbf{R}^m,\mathbf{R}^n)$.

\textbf{Theorem 3.4}\quad
$B(\mathbf{R}^m,\mathbf{R}^n)=\bigcup^n\limits_{k=0}F_k(m\geq n),
B(\mathbf{R}^m,\mathbf{R}^n)=\bigcup^n\limits_{k=0}F_k(m<n),F_k$ is
a smooth and path connected submanifold in
$B(\mathbf{R}^m,\mathbf{R}^n)$, and especially, dim$F_k=(m+n-k)k$
for $k=0,1,\cdots,n.$

\textbf{Proof}\quad We need only to prove the formula
dim$F_k=(m+n-k)k$ for $k=0,1,\cdots,n$ since otherwise  the  theorem
is  immediate from Theorem 3.3. Let $T=\{T_{i,j}\}^{m,n}_{i,j=1}$
for any $T\in B(\mathbf{R}^m,\mathbf{R}^n),I_k=\{T\in
B(\mathbf{R}^m,\mathbf{R}^n):t_{i,j}=0$ except $t_{i,i}=1, 1\leq
i\leq k\}$, and $I^{+}_k=\{(s_{i,j}\}^{m,n}_{i,j=1}\in
B(\mathbf{R}^n,\mathbf{R}^m):s_{i,j}=0$ except $s_{i,i}=1, 1\leq
i\leq k\}$.

Obviously, $I_kI^+_kI_k=I_k$ and $I^{+}_kT_kT^+_k=T^+_k.$ So
$$P^{R(I_k)}_{N(I^+_k)}=I^*_n-I_kI^+_k=\{\{s_{i,j}\}^n_1\in
B(\mathbf{R}^n):s_{i,j}=0\quad{\rm except}\quad  s_{i,i}=1,k+1\leq
i\leq n\},$$
and
$$P^{R(I^+_k)}_{N(I_k)}=I^*_m-I^+_kI_k=\{\{s_{i,j}\}^m_1\in
B(\mathbf{R}^m):s_{i,}=0\ {\rm exeept}\quad s_{i,i}=1,k+1\leq i\leq
m\},$$ where $I^*_n,I^*_m$ denote the  identity on $\mathbf{R}^n$
and $\mathbf{R}^m$, respectively. Let
$$M^+=\{P^{R(I_k)}_{N(I^+_k)}TP^{R(I^+_k)}_{N(I_k)}:\forall  T\in
B(\mathbf{R}^m,\mathbf{R}^n)\}.$$ By direct computing
$$M^+=\{\{\lambda_{i,j}\}^{n,m}_{i,j=1}:\lambda_{i,j}=0\quad{\rm
except}\quad\lambda_{i,i}=t_{i,i}, k+1\leq i\leq n\quad{\rm
and}\quad k+1\leq j\leq m\}.$$ so that dim$M^+=(n-k)(m-k).$ By Lemma
4.1 in [Ma4], $B(\mathbf{R}^m,\mathbf{R}^n)=M(I_k)\oplus M^+$ and
so, dim$M(I_k)=m\times n-(m-k)(n-k)=(m+n)-k)k.$ Due to  $F_k$ being
path connected, by Theorem 3.1, one can conclude ${\rm dim}F_k={\rm
dim}M(I_k)=(m+n-k)k.$ \quad$\Box$


Let $G(\cdot)$ denote the set of all splitling subspaces in the
Banach space in the parentheses, $U_{E}(R) = \left\{H\in G(E) : E =
R \oplus H\right\}$ for any $R\in G(E) $, and $U_{F}(S) = \{L\in
G(E) :$ $ F = S \oplus L\}$ for any $S\in G(E) $. In order to
consider of more general results then  that  of the previous
theorems 3.1 and 3.2, we introduce the equivalent relation as
follows.


{\bf Definition 3.1}\quad{\it  $T_{0}$ and $T_{*}$ in $B^{+}(E,F)$
are said to be equivalent provided there exist finite number of subspaces
$N_{1}, \cdots, N_{m}$ in $G(E)$, and $F_{1}, \cdots, F_{n}$ in $G(F)$
such that all
$$
U_{E}(N(T_{0}))\cap U_{E}(N_{1}),\cdots,U_{E}(N_{m})\cap U_{E}(N(T_{*}))
$$
and
$$
U_{F}(R(T_{0}))\cap U_{F}(F_{1}),\cdots,U_{F}(F_{n})\cap U_{F}(R(T_{*}))
$$
are non-empty. For abbreviation, write it as $T_{0}\sim T_{*},$ and let
$\widetilde{T}$ denote the equivalent class generated by $T$ in $B(E,F).$}


{\bf Theorem 3.5}\quad{\it  $\widetilde{T_{0}}$ for any $T_{0}\in
B^{+}(E,F)$ with $\mathrm{codim } R(T_{0}) > 0$ is path connected.}

{\bf Proof}\quad Assume that $T_{*}$ is any operator in $\widetilde{T_{0}}$, and
$$
R_{1}\in U_{E}(N(T_{0}))\cap U_{E}(N_{1}),
\cdots,R_{m}\in U_{E}(N_{m-1})\cap U_{E}(N_{m}),
       R_{m+1}\in U_{E}(N_{m})\cap U_{E}(N_(T_{*})).\eqno(3.7)
$$

We define by induction:
$$
T_{k} = T_{k-1}P^{N_{k}}_{R_{k}},\,\,k=1,2,\cdots,m.
$$
It is easy to observe
$$
R(T_{k}) = R(T_{0})\,\mathrm{and }\, N(T_{k}) = N_{k}, k=1,\cdots,m.\eqno(3.8)
$$
Indeed, $R(T_1)=T_0R(P^{N_1}_{R_1})=T_0R_1=R(T_0)$ and$N(T_1)=\{x\in
E:P^{N_1}_{R_1}x\in N(T_0)\}=\{x\in E:P^{N_1}_{R_1}x=0\}=N_1$
because of $E=R_1\oplus N_1$; similarly, by induction one can
conclude $R(T_k)=R(T_0)$ and $N(T_k)=N_k$ for $k=1,2,\cdots,m.$ So
$T_k\in\tilde{T}_0$ for $k=1,2,\cdots,m.$

Let $S_k=\{T\in G_d(R_k):R(T)=R(T_0)\}=\{T\in B(E,F):E=N(T)\oplus
R_k$ and $R(T)=R(T_0)\}$ for $k=1,2,\cdots,m.$ Evidently $S_k\subset
\tilde{T}_0, k=1,2,\cdots,m.$ In fact, by (3.7) $R_1\in
U_E(N(T_0))\cap U_E(N_1),\cdots,R_k\in U_E(N_{k-1})\cap U_E(N_k);$
while $R_k\in U_E(N_k)\cap U_E(N(T))$ and $R(T)=R(T_0);$ so that
$S_k\subset\tilde{T}_0,k=1,2,\cdots,m.$ Next go to show that $T_0$
and $T_m$ are path connected in $\tilde{T}_0$. Due to $N_{k-1}\oplus
R_k)=N(T_{k-1})\oplus R_k=N_k\oplus R_k.$ Theorem 2.4 shows that
$T_{k+}$ and $T_k$ are path connected in $S_k$, and so in
$\tilde{T}_0$ for $k=1,2,\cdots,m$. Therefore $T_0$ and $T_m$ are
path connected in $\tilde{T}_0$.

 Let $M_k=\{T\in C_d(R_k):\mathbf{R}(T)=R(T_k)\}.$
By (3.8)
$$M_k=R(T_k)\}=\{T\in
C_d(R_k):R(T)=R(T_0)\},\quad k=1,2,\cdots,m,$$ so that
$M_k\subset\tilde{T}_0$. Note $E=N(T_{k-1})\oplus R_k=N_k\oplus
R_k.$ Then by Theorem 2.4, $T_k=T_{k-1}P^{N_k}_{R_k}$ and $T_{k-1}$
are path connected in $M_{k-1}=\{T\in
C_d(R_k):R(T)=R(T_{k-1})\}\subset\tilde{T}_0, k=1,2,\cdots,m.$ This
shows that $T_0$ and $T_m$ are poth connected in $\tilde{T}_0$.
Finally go to prove that $T_m$ and $T_*$ are path connected in
$\tilde{T}_0.$

In other hand, assume
$$
S_{1}\in U_{F}(R(T_{0})) \cap U_{F}(F_{1}),\cdots,
S_{n}\in U_{F}(F_{n-1}) \cap U_{F}(F_{n}),
S_{n+1}\in U_{F}(F_{n}) \cap U_{F}(R(T_{*})).\eqno(3.9)
$$

 Write $T_{m,0}=T_m.$
We define  by induction,
$$
T_{m,i}=P_{F_{i}}^{s_{i}}T_{m,i-1},\, i=1, 2,\cdots,n.
$$

According to the equality  $R(T_{m})= R(T_{0})$ in (3.8), we infer
$$
R(T_{m,i})=F_{i}\,\mathrm{and}\,N(T_{m,i})=N(T_{m}),\,\,i=1,\cdots,n.\eqno(3.10)
$$
In fact wsmite $F_0=\mathbf{R}(T_m)$, then  by (3.9)
$$
F = F_{i-1}\oplus S_{i} =  F_{i} \oplus S_{i}, \,i = 1, 2,\cdots,n;
$$
and so
$$N(T_{m,i})=N(T_m)\quad{\rm and}\quad
R(T_{m,i})=F_i,i=1,2,\cdots,n.$$

Thus, take $T_{m,i-1},F_i$ and $S_i$ for $i=1,2,n$ in the places of
$T_0,F_*$ and $N$, in Theorem 2.3, respectively, then the theorem
shows that  $T_{m,k}$ and $T_{m,k-1}$ are path connected in
$\left\{T\in C_{r}(S_{k}): N(T) = N_{m}\right\} \subset
\widetilde{T_{0}}, k=1,2,\cdots,n,$ so that $T_{m}$ and $T_{m,n}$
are path connected in $\widetilde{T_{0}}$. Since  $T_{0}$ and
$T_{m}$ are path connected in $\widetilde{T_{0}},$  $T_{0}$ and
$T_{m,n}$ are path connected in $\widetilde{T_{0}},$ so  the proof
of the theorem  reduces to that of $T_{m,n}$ being path connected
with $T_{*}$ in $\widetilde{T_{0}},$ where $T_{m,n}$ and $T_{*}$
satisfy
$$
E = N(T_{m,n})\oplus R_{m+1}=  N(T_{*})\oplus R_{m+1} \,\, \mathrm{because\, of }\,\,
N(T_{m,n}) = N_{m},
$$
and
$$
F = R(T_{m,n})\oplus S_{n+1}=  R(T_{*})\oplus S_{n+1} \,\, \mathrm{because\, of }\,\,
R(T_{m,n}) = F_{n}.
$$
$$\eqno(3.11)$$
The first equality in (3.11) shows $T_{*}\in C_{\alpha}(R_{m+1})$
and $E = N_{m}\oplus R_{m+1}.$ Then by Theorem 2.4,
$T_{*}P^{N_{m}}_{R_{m+1}}$ and $T_{*} $ are path connected in
$\left\{T\in C_{\alpha}(R_{m+1}): R(T) = R(T_{*})\right\} \subset
\widetilde{T_{0}}.$

The second equality in (3.11) shows $T_{*}\in C_{r}(S_{n+1})$ an d
$F = F_{n}\oplus S_{n+1}.$ Then by Theorem 2.3,
$P^{S_{n+1}}_{F_{n}}T_{*}$ and $T_{*} $ are path connected in
$\left\{T\in C_{r}(S_{n+1}): N(T) = N(T_{*})\right\} \subset
\widetilde{T_{0}}.$ Combining the preceding results, we conclude
that $P^{S_{n+1}}_{F_{n}}T_{*}P^{N_{m}}_{R_{n+1}}$ and $T_{*}$ are
path connected in $\widetilde{T_{0}}.$  So far, the proof of the
theorem reduces to that of
 $P^{S_{n+1}}_{F_{n}}T_{*}P^{N_{m}}_{R_{m+1}}$ being path connected with $T_{m,n}$ in $\widetilde{T_{0}}.$

According to (3.10) we have
$$
E = N_{m}\oplus R_{m+1}\,\,\mathrm{and}\,\,F = F_{n}\oplus S_{n+1}
$$
where $N_{m}=N(T_{m,n})$ and $R(T_{m,n})=F_{n}.$

Let
$$T_{m,n}^{+}y =\left\{\begin{array}{rcl}\left(T_{m,n}|_{R_{m+1}}\right)^{-1}y, & & {y\in F_{n}},\\
0, & & {y\in S_{n+1}}.
\end{array}\right.$$

Obviously,
$$
P^{S_{n+1}}_{F_n}T_*P^{N_m}_{R_{m+1}}=P^{S_{n+1}}_{F_{n}}T_{*}T^{+}_{m,n}T_{m,n}.\eqno(3.12)
$$
We claim
$$
P^{S_{n+1}}_{F_{n}}T_{*}T^{+}_{m,n}|_{F_{n}}\in B^{X}(F_{n})
$$
$ P^{S_{n+1}}_{F_{n}}T_{*}T^{+}_{m,n}y=0 \,\,\mathrm{for}\,\,y\in
F_{n}\Leftrightarrow T_{*}{T^{+}_{m,n}}y=0 $ because of (3.11)
$\Leftrightarrow T^{+}_{m,n}y=0 $ since $T^{+}_{m,n}y\in R_{m+1}
\Leftrightarrow y=0;$ while from  the assumption in (3.7)and (3.8),
$R_{m+1}\in  U_{E}(N_{m}) \cap U_{E}(N(T_{*})) $ and $S_{n+1} \in
U_{F}(F_{n}) \cap U_{F}(R(T_{*})),$ it follows that  there exist
$r\in R_{m+1}$ and $y_{0}\in F_{n},$ for any $y\in F_{n}$ such that
$ P^{S_{n+1}}_{F_{n}}T_{*}r =
y\,\,\mathrm{and}\,\,T^{+}_{m,n}y_{0}=r,\,i.e.,
P^{S_{n+1}}_{F_{n}}T_{*}T^{+}_{m,n} $ is surjective.

It is well known that $P^{S_{n+1}}_{F_{n}}T_{*}T^{+}_{m,n}|_{F_{n}}$
is path connected with some one of  $I_{F_{n}}$ and  $(-I_{F_{n}})$
in $B^{X}(F_{n}).$ Hence
$P^{S_{n+1}}_{F_{n}}T_{*}T^{+}_{m,n}P^{S_{n+1}}_{F_{n}}$ is path
connected with some one  $P^{S_{n+1}}_{F_{n}}$ and
$(-P^{S_{n+1}}_{F_{n}})$ in the set $S=$
 $\left\{F\in B(F_{n}): R(T) =
F_{n}\,\, \mathrm{and}\,\, N(T) = S_{n+1}\right\}$ $\subset
\left\{T\in C_{r}(S_{n+1}): N(T)=S_{n+1}\right\}.$ If
$P^{S_{n+1}}_{F_{n}}T_{*}T^{+}_{m,n}P^{S_{n+1}}_{F}$ is path
connected with $(-P^{S_{n+1}}_{F_{n}})$, then by Theorem 2.2, it is
path connected with $P^{S_{n+1}}_{F_{n}}$ in $S$. Thus
$P^{S_{n+1}}_{F_{n}}T_{*}T^{+}_{m,n}P^{S_{n+1}}_{F_{n}}$ is also
path connected with $P^{S_{n+1}}_{F_{n}}$ in $\left\{T\in
C_{r}(S_{n+1}): N(T)= S_{n+1}\right\}$.
 Note the equality (3.10). Finally, by the simular way to that in
 the proof of Theorem 3.2, one can infer that
 $P^{S_{n+1}}_{F_n}T_*P^{N_m}_{R_{m+1}}$ and $T_m$ are path connected
 in $\tilde{T}_0.$ The proof ends.\quad$\Box$


It is easy to see that $\widetilde{T} = F_{k}$ for any $T\in F_{k}
(k< \infty)$  as well as
 $\widetilde{T} = \Phi_{m,n}$ for any $T\in \Phi_{m,n} ,\, n
> 0.$

\vskip 0.1cm
\begin{center}{\bf References}
\end{center}
\vskip -0.1cm
\medskip
{\footnotesize
\def\REF#1{\par\hangindent\parindent\indent\llap{#1\enspace}\ignorespaces}

\REF{[Abr]}\ R. Abraham, J. E. Marsden, and T. Ratin, Manifolds,
tensor analysis and applications, 2nd ed., Applied Mathematical
Sciences 75, Springer, New York, 1988.

\REF{[An]}\ V. I. Arnol'd, Geometrical methods in the theory of
ordinary differential equations, 2nd ed., Grundlehren der
Mathematischen Wissenschaften 250, Springer, New York, 1988.

\REF{[Bo]}\ B. Booss,D.D.Bleecker, Topology and Analysis: the
Atyah-Singer Index Formula and Gauge-Theoretic physics, New York:
Springer-Verlag 1985.

\REF{[Caf]}\ V. Cafagra, Global invertibility and finite
solvability, pp. 1-30 in Nonlinear functional analysis (New York,
NJ, 1987), edited by P. S. Milojevic, Lecture Notes in Pure and
Appl. Math. 121, Dekker, New York, 1990.

\REF{[Ma1]}\ Jipu Ma, Dimensions of Subspaces in a Hilbert Space and
Index of Semi-Fredholm Operators, Sci, China, Ser. A 39:12(1986).

\REF{[Ma2]}\ Jipu Ma, Generalized Indices of Operators in B(H), Sci,
China, Ser. A , 40:12(1987).

\REF{[Ma3]}\ Jipu Ma, Complete Rank Theorem in Advanced Calculus and
Frobenius Theorem in Banach Space, arxiv:1407.5198v5[math.FA]23 Jan.
2015.

\REF{[Ma4]}\ Jipu Ma Frobenius Theorem in Banach Space,  submited to
Ivent. Math..(Also see) arXiv: submit/1512157[math.FA] 19 Mar 2016).

  1. Department of Mathematics, Nanjing University, Nanjing, 210093,
  P. R. China

  2. Tseng Yuanrong Functional Analysis Research Center, Harbin Normal
  University, Harbin, 150080, P. R. China

 E-mail address: jipuma@126.com; jipuma1935@sina.com
}
\end{document}